\documentclass[10pt]{article}
\usepackage{amsfonts}
\usepackage{amsmath}
\usepackage{amssymb}
\usepackage{graphicx}%

\newtheorem{theorem}{Theorem}

\newtheorem{corollary}[theorem]{Corollary}

\newtheorem{lemma}[theorem]{Lemma}
\newtheorem{notation}[theorem]{Notation}

\newtheorem{proposition}[theorem]{Proposition}
\newtheorem{remark}[theorem]{Remark}

\def \NN {{\mathbb N}}
\newcommand{\bpr}{\noindent\textbf{Proof\/. }}
\newcommand{\epr}{\hspace*{\fill}$\Box$\medskip}
\begin{document}

\title{The $\{-3\}$-reconstruction and the $\{-3\}$-self duality of tournaments}
\maketitle


\centerline{Mouna Achour$^{(a)}$, Youssef Boudabbous$^{(a)}$,
Abderrahim Boussa\"{\i}ri$^{(b)}$}

\vskip0.5cm (a) D\'epartement de Math\'ematiques, Facult\'e des
Sciences de Sfax, Universit\'e de Sfax, BP $1171$, $3038$ Sfax,
Tunisie. Fax: (00216) 74.27.44.37\newline E-mail: mouna\underline{
}achour@yahoo.fr, youssef\underline{ }boudabbous@yahoo.fr
\vskip0.001 cm (b) Facult\'e des Sciences A\"{\i}n Chock,
D\'epartement de Math\'ematiques et Informatique, Km 8 route d'El
Jadida, BP 5366 Maarif, Casablanca, Maroc. E-mail:
aboussairi@hotmail.com

\begin{abstract}
Let $T = (V,A)$ be a (finite) tournament and $k$ be a non negative
integer. For every subset $X$ of $V$ is associated the subtournament
$T[X] = (X,A\cap (X \times X))$ of $T$, induced by $X$. The dual
tournament of $T$, denoted by $T^{\ast}$, is the tournament obtained
from $T$ by reversing all its arcs. The tournament $T$ is self dual
if it is isomorphic to its dual. $T$ is $\{-k\}$-self dual if for
each set $X $ of $k$ vertices, $T[V\setminus X]$ is self dual. $T$
is strongly self dual if each of its induced subtournaments is self
dual. A subset $I$ of $V$ is an interval of $T$  if for  $a, b \in I
$ and for   $x\in V\setminus I$, $ (a,x) \in A$ if and only if $
(b,x) \in A$. For instance, $\emptyset$, $V$ and $\{x\}$, where
$x\in V$, are intervals of $T$ called trivial intervals. $T$ is
indecomposable if all its intervals are trivial; otherwise, it is
decomposable. A tournament $T'$, on the set $V$, is
$\{-k\}$-hypomorphic to $T$ if for each set $X $ on $k$ vertices,
$T[V\setminus X]$ and $T'[V\setminus X]$ are isomorphic. The
tournament $T$ is $\{-k\}$-reconstructible if each tournament
$\{-k\}$-hypomorphic to $T$ is isomorphic to it.\\
Suppose that $T$ is decomposable and $\mid V\mid \geq 9$. In this
paper, we begin by proving the equivalence between the $\{-3\}$-self
duality and the strong self duality of $T$. Then we characterize
each tournament $\{-3\}$-hypomorphic to $T$. As a consequence of
this characterization, we prove that if there is no interval $X$ of
$T$ such that $T[X]$ is indecomposable and $\mid V \setminus X\mid
\leq 2$, then $T$ is $\{-3\}$-reconstructible. Finally, we conclude
by reducing the $\{-3\}$-reconstruction problem to the
indecomposable case (between a tournament and its dual). In
particular, we find and improve, in a less complicated way, the
results of \cite{Boud et Bouss} found by Y. Boudabbous and A.
Boussa\"{\i}ri.
\end{abstract}

\section{Introduction}

\subsection{Preliminaries on tournaments} A (finite) \emph{tournament} $T$ consists of
a finite set $V$ of \emph{vertices} with a prescribed collection $A$
of ordered pairs of distinct vertices, called the set of \emph{arcs}
of $T$, which satisfies: for $x,y\in V$ with $x\neq y$, $(x,y)\in A$
if and only if $(y,x)\not \in A$. Such a tournament is denoted by
$(V,A)$. If $(x,y)$ is an arc of $T$, then we say that $x$
\emph{dominates} $y$ (symbolically $x\rightarrow y$). The
\emph{dual} of the tournament $T$ is the tournament
$T^{\ast}=(V,A^{\ast})$ defined by: for all $x$, $y\in V$, $(y,x)\in
A^{\ast}$ if and only if $(x,y)\in A$. The tournament $T$ is
\emph{transitive} or a \emph{linear order} provided that for any
$x,y,z\in V$, if $(x,y)\in A$ and $(y,z)\in A$, then $(x,z)\in A$.
For example, a total order on a finite set $E$ can be identified to
a transitive tournament with a vertex set $E$ in the following way:
for $x,y\in E$ with $x\neq y$, $x\rightarrow y$ if and only if
$x<y$. The tournament corresponding to the usual order on
$\{1,\ldots,n\}$ (where $n \in \NN^{\ast}$) is denoted by $O_{n}$.
An \emph{almost transitive tournament} is a tournament obtained from
a transitive tournament with at least three vertices by reversing
the arc formed by its two extremal vertices.

For every finite sets $E$ and $F$, we denote $E \subset F$ when $E$
is a subset of $F$ and $\mid E\mid$ the cardinality of $E$.

Given a tournament $T=(V,A)$, for each subset $X$ of $V$ we
associate the \emph{subtournament} of $T$ induced by $X$, that is
the tournament $T[X] = (X,A\cap\\(X\times X))$. For convenience, the
subtournament $T[V\backslash X]$ is denoted by $T-X$, and by $T-x$
whenever $X=\{x\}$.

Let $T=(V,A)$ be a tournament, a subset $I$ of $V$ is an
\emph{interval} of $T$ if for every $x\in V\setminus I$, $x$
dominates or is dominated by all elements of $I$. For instance,
$\emptyset$, $V$ and $\{x\}$, where $x\in V$, are intervals of $T$
called \emph{trivial} intervals. A tournament is
\emph{indecomposable} if all its intervals are trivial; otherwise,
it is
\emph{decomposable}. For example, the tournament $C_{3}%
=(\{1,2,3\},\{(1,2),(2,3),(3,1)\})$ is indecomposable, whereas, the
tournaments $C_{4}=(\{1,2,3,4\},\{(1,2),(2,3),(3,4),(4,1),(3,1),(2,4)\})$,
$\delta^{+}=(\{1,2,3,4\},\\ \{(1,2),(2,3),(3,1),(1,4),(2,4),(3,4)\})$
and $\delta^{-}=(\delta^{+})^{\ast}$ are decomposable.

Given two tournaments $T=(V,A)$ and
$T^{\prime}=(V^{\prime},A^{\prime})$, an \emph{isomorphism} from $T$
onto $T^{\prime}$ is a bijection $f$ from $V$ onto $V^{\prime}$
satisfying: for any $x,y\in V$, $(x,y)\in A$ if and only if
$(f(x),f(y))\in A^{\prime}$. The tournaments $T$ and $T^{\prime}$
are \emph{isomorphic} if there exists an isomorphism from one onto
the other. This is denoted by $T\sim T^{\prime}$. A tournament
$T^{\prime}$ \emph{embeds} into a tournament $T$ (or $T$
\emph{embeds} $T^{\prime}$ ), if $T^{\prime}$ is isomorphic to a
subtournament of $T$. A \emph{$3$-cycle} (resp. \emph{$4$-cycle}) is
a tournament which is isomorphic to $C_{3}$ (resp. $C_{4}$).
Moreover, a \emph{positive diamond} (resp. \emph{negative diamond})
is a tournament that is isomorphic to $\delta^{+}$ (resp.
$\delta^{-}$). A \emph{diamond} is a positive or a negative diamond.
For convenience, a set $X$ of vertices of a tournament $T$ is called
diamond of $T$ if $T[X]$ is a diamond.

\subsection{Self duality and reconstruction} A tournament $T$ on a set $V$ is \emph{self dual} if $T$ and
$T^{\ast}$ are isomorphic, it's \emph{strongly self dual} if for
every subset $X$ of $V$, $T[X]$ and $T^{\ast}[X]$ are isomorphic.
For each non negative integer $k$, the tournament $T$ is
\emph{$(\leq k)$-self dual} whenever for every set $X$ of at most
$k$ vertices, the subtournament $T[X]$ is self dual. It is easy to
see that a transitive tournament or an almost transitive tournament
is strongly self dual. Conversely, Reid and Thomassen
\textrm{\cite{Reid et Thom}} was proved that a strongly self dual
tournament with at least $8$ vertices is transitive or almost
transitive. This result was used by K. B. Reid and C. Thomassen
\cite{Reid et Thom} in order to characterize the pair of
\emph{hereditarily isomorphic} tournaments, that is, the pair of
tournaments $T$ , $T^{\prime}$ on a set $V$ such that for every
subset $X$ of $V$, the subtournaments $T[X]$ and $T^{\prime}[X]$ are
isomorphic. A relaxed version of this notion is the following.
Consider two tournaments $T$ and $T^{\prime}$ on the same vertex set
$V$, with $\mid V\mid=n\geq2$. Let $k$ be an non negative integer $k
$ with $k \leq n$. The tournaments $T$ and $T^{\prime}$ are
\emph{$\{k\}$-hypomorphic}, whenever for every set $X$ of $k$
vertices, the subtournaments $T[X]$ and $T^{\prime}[X]$ are
isomorphic. If $T$ and $T^{\prime}$ are $\{n-k\}$-hypomorphic, we
say that $T$ and $T^{\prime}$ are \emph{$\{-k\}$-hypomorphic}. Let
$F$ be a set of integers. The tournaments $T$ and $T^{\prime}$ are \emph{$F$%
-hypomorphic}, if for every $p\in F$, $T$ and $T^{\prime}$ are $\{p\}$%
-hypomorphic, in particular, if $F=\{0,\ldots,k\}$, we say that $T$
and $T^{\prime}$ are \emph{$(\leq k)$-hypomorphic}. For example,
every two tournaments on the same vertex set with at least $2$
elements are $(\leq2)$-hypomorphic. The tournament $T$ is
\emph{$F$-reconstructible} provided that every tournament
$F$-hypomorphic to $T$ is isomorphic to $T$. This notion was
introduced by R. Fra\"{\i}ss\'{e} \cite{Fraiss} in 1970.

In 1972, G. Lopez (\textrm{\cite{Lopez 2res},\cite{lopez suite}})
showed that a tournament, with at least $6$ vertices, is
$(\leq6)$-reconstructible (see also \textrm{\cite{Lopez indéf}}). It
follows from a "Combinatorial Lemma" of M. Pouzet (see Section $2$)
that a tournament, with at least $12$ vertices, is
$\{-6\}$-reconstructible.

On the other hand, P. K. Stockmeyer \cite{Stock} showed that the
tournaments are not, in general, $\{-1\}$-reconstructible,
invalidating so the conjecture of Ulam \cite{Ulam} for tournaments.
Then, M. Pouzet (\cite{Bondy},\cite{Bondy et Hemm}) proposed the
$\{-k\}$-reconstruction problem of tournaments. P. Ille
\textrm{\cite{Ille}} established that a tournament with at least
$11$ vertices is $\{-5\}$-reconstructible. G. Lopez and C. Rauzy
(\textrm{\cite{Lopez+Rauzy I},\cite{Lopez+Rauzy II}}) showed that a
tournament with at least $10$ vertices is $\{-4\}$-reconstructible.
The $\{-k\}$-reconstruction problem of tournaments is still open for
$k\in \{2,3\}$.

In 1995, Y. Boudabbous and A. Boussa\"{\i}ri  \cite{Boud et Bouss} studied the $\{-3\}$
-recons-truction of decomposable tournaments, for which they give a
partial positive answer.

In this paper, we find and improve, in a less complicated way, the
results of this study. In order to present our results, we need the
Gallai's decomposition theorem \textrm{\cite{Gallai} }.
\subsection{Gallai's decomposition} Given a tournament $T=(V,A)$, we define on $V$ a
binary relation $\mathcal{R}$ as follows: for all $x\in V$,
$x\mathcal{R}x$ and for $x\neq y\in V$, $x\mathcal{R}y$ if there
exist two integers $n,\, m \geq 1$ and two sequences
$x_{0}=x,\ldots,x_{n}=y$ \ and $y_{0}=y,\ldots,y_{m}=x$ of vertices
of $T$ such that $x_{i}\longrightarrow x_{i+1}$ for $i=0,\ldots,n -
1$ and $y_{j}\longrightarrow y_{j+1}$ for $j=0,\ldots,m - 1$.
Clearly, $\mathcal{R}$ is an equivalence relation on $V$. The
equivalence classes of $\mathcal{R}$ are called the \emph{strongly
connected components} of $T$. A tournament is then \emph{strongly
connected} if it has at most one strongly connected component,
otherwise, it is \emph{non-strongly connected}.

The next result is due to J. W. Moon \cite{Moon}.

\begin{lemma}
\label{Moon}\textrm{\cite{Moon}} Given a strongly connected
tournament $T=(V,A)$ with $n\geq3$ vertices, for every integer
$k\in\{3,\ldots,n\}$ and for every $x\in V$, there exists a subset
$X$ of $V$ such that $x\in X$, $\mid X\mid=k$ and the subtournament
$T[X]$ is strongly connected.
\end{lemma}

Let $T$ be a tournament on a set $V$. A partition $\mathcal{P}$
of $V$ is an \emph{interval partition} of $T$ if all the elements of
$\mathcal{P}$ are intervals of $T$. It ensues that the elements of
$\mathcal{P}$ may be considered as the vertices of a new tournament,
the \emph{quotient} $T/\mathcal{P}=(\mathcal{P},A/\mathcal{P})$
\emph{of $T$ by $\mathcal{P}$}, defined in the following way: for
any $X\neq Y\in\mathcal{P}$, $(X,Y)\in A/\mathcal{P}$ if $(x,y)\in
A$, for $x\in X$ and $y\in Y$. On another hand, a subset $X$ of $V$
is a \emph{strong interval} of $T$ provided that $X$ is an interval
of $T$ and for every interval $Y$ of $T$, if $X\cap Y\neq\emptyset$,
then $X\subset Y$ or $Y\subset X$. Here, for each tournament
$T=(V,A)$ with $\mid V\mid\geq2$, $\mathcal{P}(T)$ denotes the
family of maximal, strong intervals of $T$, under the inclusion,
amongst the strong intervals of $T$ distinct from $V$. Clearly,
$\mathcal{P}(T)$ realizes an interval partition of $T$.

Consider a tournament $H = (V,A)$. For every $x \in V$ is associated
the tournament $T_{x} = (V_{x},A_{x})$ such that the $V_{x}$'s are
mutually disjoint. The lexicographical sum of $T_{x}$'s over $H$ is
the tournament $T$ denoted by $H(T_{x};x\in V)$ and defined on the
union of $V_{x}$'s as follows: given $u \in V_{x}$ and $v \in
V_{y}$, where $x$, $y\in V$, $(u,v)$ is an arc of $T$ if either $x =
y$ and $(u,v)\in A_{x}$ or $x \neq y$ and $(x,y)\in A$. This
operation consists in fact to replace every vertex $x$ of $V$ by
$T_{x}$ so that $V_{x}$ becomes an interval; we say that the vertex
$x$ is \emph{dilated} by $T_{x}$. For example, an almost transitive
tournament is obtained from a $3$-cycle by dilating one of its
vertices by a transitive tournament.

The Gallai's decomposition theorem \textrm{\cite{Gallai} }consists
in the following examination of the quotient $T/\mathcal{P}(T)$.

\begin{theorem} (\textrm{\cite{BILT},\cite{Gallai}})
\label{gallai theo} Let $T$ be a tournament with at least two
vertices.

\begin{enumerate}
\item The tournament $T$ is non-strongly connected if and only if
$T/\mathcal{P}(T)$ is transitive. In addition, if $T$ is non-strongly connected, then $\mathcal{P}%
(T)$ is the family of the strongly connected components of $T$.

\item The tournament $T$ is strongly connected if and only if $T/\mathcal{P}%
(T)$ is indecomposable and $\mid\mathcal{P}(T)\mid\geq3$.
\end{enumerate}
\end{theorem}

\bigskip

We complete this subsection by the following notation.

\begin{notation}
For every tournament $T$ defined on a vertex set $V$ with at least
two elements, we associate the partition
$\widetilde{\mathcal{P}}(T)$ of $V$ defined from $\mathcal{P}(T)$ as
follows:

\begin{itemize}
\item If $T$ is strongly connected,
$\widetilde{\mathcal{P}}(T)=\mathcal{P}(T)$.

\item If $T$ is non-strongly connected, a subset $A$ of $V$ belongs to
$\widetilde{\mathcal{P}}(T)$ if and only if either
$A\in\mathcal{P}(T)$ and $\mid A\mid\geq2$, or $A$ is a maximal
union of consecutive vertices of the transitive tournament
$T/\mathcal{P}(T)$ which are singletons.
\end{itemize}
\end{notation}

\subsection{Statement of the results}

In this paper, we begin by proving the following theorem (Section
$3$). That improves the result obtained by K. B. Reid and C.
Thomassen \cite{Reid et Thom}. (Note that, the study of \cite{Reid
et Thom} is easly reduced to the decomposable case).

\begin{theorem}
\label{theorem principal moins 3 autodual} A decomposable tournament which has
at least $9$ vertices is $\{-3\}$-self dual if and only if it is either transitive or
almost transitive.
\end{theorem}

Then, we characterize each tournament $\{-3\}$-hypomorphic to a
decomposable tournament with at least $9$ vertices (Section $4$).
For the statement, we need additional notation.

\begin{notation}
Consider a set $P$ of non zero integers, an integer $n\geq 6$ and \\
$q\in \{2,3\}$. We denote by:

\begin{itemize}
\item $I_{n,P}$, the class of tournaments with $n$ vertices which are
indecomposable, not self dual and $\{p\}$-self dual for every $p \in
P$; these tournaments being considered up to an isomorphism.

\item $C_{3}(I_{n,P})$ (resp. $O_{q}(I_{n,P})$), the class of tournaments
with $n+2$ (resp. $n+q-1$) vertices obtained, from the $3$-cycle
$C_{3}$ (resp. from the transitive tournament $O_{q}$), by dilating
one of its vertices by a tournament belonging to the class
$I_{n,P}$; these tournaments being considered up to an isomorphism.

\item For each integer $m\geq 8$, we denote by $\Omega_{m}$ the union
$C_{3}(I_{m - 2,\{-1,-2,-3\}})\\ \cup O_{3}(I_{m - 2,\{-1,-2,-3\}}) \cup O_{2}(I_{m - 1,\{-2,-3\}})$.
\end{itemize}
\end{notation}

 Here is the characterization.

\begin{theorem}
\label{Theo moins 3 reconstruction}Consider a decomposable tournament $T$
with $n\geq9$ vertices and let $T'$ be a tournament $\{-3\}$-hypomorphic to $T$. Then, we have:
\begin{enumerate}
  \item $\widetilde{\mathcal{P}}(T^{\prime})=\widetilde
{\mathcal{P}}(T)$ and one of the following situations is achieved.
\begin{description}
  \item[$(a)$] $T$ is almost transitive and $T' \sim T$.
  \item[$(b)$] $T$ is not almost transitive,  $T/\widetilde{\mathcal{P}}(T)=T'/\widetilde
{\mathcal{P}}(T)$ and one of the following situations is achieved.
\begin{description}
  \item[$(i)$] $T \not\in \Omega_{n}$, for every $X\in \widetilde{\mathcal{P}}(T)$, $T'[X]\sim T[X]$ and $T'\sim T$.
  \item[$(ii)$] $T\in \Omega_{n}$ and $T' \not\sim T$.
\end{description}
\end{description}
  \item $T$ is not $\{-3\}$-reconstructible if and only if $T \in \Omega_{n}$.
\end{enumerate}
\end{theorem}

We deduce the following results.

\begin{corollary}
\label{ajoutth1}Let $T$ be a decomposable tournament on a set $V$,
with $\mid V\mid \geq9$. If there is no interval $X$ of $T$ such
that $T[X]$ is indecomposable and $\mid V\setminus X\mid\leq2$, then
$T$ is $\{-3\}$-reconstructible.
\end{corollary}

\begin{corollary}
\label{theo moin3 recon} If there exists an integer $n_{0}\geq7$ for which the
indecomposable tournaments with at least $n_{0}$ vertices are $\{-3\}$-reconstructible,
then the tournaments with at least $n_{0} + 2$ vertices are
$\{-3\}$-reconstructible.
\end{corollary}

By Theorem \ref{theo inversion} (Section $2$) and Corollary
\ref{theo moin3 recon}, we reduce the $\{-3\}$-reconstruction
 problem of tournaments to the indecomposable case between a tournament and its dual.

\section{ Preliminary results}

In this section, we recall and prove some results which will be used
in next sections.

First, concerning the decomposability, we recall the following
notation, lemma and corollary.

\begin{notation}
Given a tournament $T=(V,A)$, for each subset $X$ of $V$, such that
$\mid\!X\!\mid\geq3$ and $T[X]$ is indecomposable, we associate the
following subsets of $V\setminus X$.

\begin{itemize}
\item $Ext(X)$ is the set of $x\in V\setminus X$ such that $T[X\cup\{x\}]$ is indecomposable.

\item $[X]$ \ is the set of $\ x\in V\setminus X\,$\ such that\ $X$ is an interval of
$T[X\cup\{x\}]$.

\item For every $u\in X$, $X(u)$ is the set of \ $x\in V\setminus X\ $such
that\ $\{u,x\}$ is an interval of $T[X\cup\{x\}]$.
\end{itemize}
\end{notation}

\begin{lemma}
\label{induicedpartition}\textrm{\cite{eur_rosen}} Let $T=(V,A)$ be
a tournament and let $X$ be a subset of $V$ such that
$\mid\!X\!\mid\geq3$ and $T[X]$ is indecomposable. The family
$Ext(X),[X]$ and $X(u)$ where $u\in X$ constitutes a partition of
$V\setminus X$. (Some elements of this family can be empty).
\end{lemma}

The next result follows from Lemma \ref{induicedpartition}.

\begin{corollary}
\label{indec}\textrm{\cite{eur_rosen}} Let $T=(V,A)$ be an
indecomposable tournament. If $X$ is a subset of $V$, such that
$\mid\!X\!\mid\geq3$, $\mid\!V\setminus X\!\mid\geq2$ and $T[X]$ is
indecomposable, then there are distinct $x, y \in V\setminus X$ such
that $T[X \cup\{x,y\}]$ is indecomposable.
\end{corollary}

Second, we consider the following remark and notation.
\begin{remark}
\label{tournaments with 4 vertices}

\begin{description}
\item[i)] Up to an isomorphism, there exist four tournaments with four
vertices: $O_{4}$, $C_{4}$, $\delta^{+}$ and $\delta^{-}$. In
addition, both of them is decomposable.

\item[ii)] Given two $\{3\}$-hypomorphic tournaments $T$ and $T^{\prime}$ with
the same vertex set $V$ with $\mid V\mid \geq 4$, $T$ and
$T^{\prime}$ are $(\leq4)$-hypomorphic if and only if for every
subset $X$ of $V$, if $T[X]$ or $T^{\prime}[X]$ is a diamond, then
$T^{\prime}[X]\sim T[X]$.

\item[iii)] Consider two $\{3\}$-hypomorphic tournaments $T$ and $T^{\prime}%
$. If $T$ is without diamonds then $T^{\prime}$ and $T$ are
$\{4\}$-hypomorphic.
\end{description}
\end{remark}

\begin{notation}
Given a tournament $T=(V,A)$, a subset $F$ of $V$ and a tournament
$H$, we denote by $S(T,H;F) = \{X\subset V;\,F\subset X$ and
$T[X]\sim H\}$ and $n(T,H;F) = \mid S(T,H;F) \mid$.
\end{notation}

Then, we recall the following "Combinatorial Lemma" of M. Pouzet
\cite{Pouz prop comb}.

\begin{lemma}
\label{lemme combinaoire}\textrm{\cite{Pouz prop comb}} Let $p$, $r$
be two positive integers, $E$ be a set of at least $p + r$ elements
and $U$, $U^{\prime}$ be two sets of subsets of $p$ elements of $E$.
If for each subset $Q$ of $E$ with  $\mid Q \mid = p + r$, the
number of the elements of $U$ which are contained in $Q$ is equal to
the number of the elements of $U^{\prime}$ which are contained in
$Q$, then for every finite subsets $P^{\prime}$ and $Q^{\prime}$ of
$E$, such that $P^{\prime}$ is contained in $Q^{\prime}$ and
$Q^{\prime}\setminus P^{\prime}$ has at least $p + r$ elements, the
number of elements of $U$ containing $P^{\prime}$ and included in
$Q^{\prime}$ is equal to the number of elements of $U^{\prime}$
containing $P^{\prime}$ and included in $Q^{\prime}$. In particular,
if $E$ has at least $2p + r$ elements, then $U$ and $U^{\prime}$ are
equal.
\end{lemma}

From the Combinatorial Lemma, we have the next corollaries.

\begin{corollary}
\label{corollaire 2 pouz} \textrm{\cite{Pouz rel non reconst}}
Consider positive integers $n$, $p$, $h$ such that $ p < n$ and
$h\leq n-p$, a tournament $H$ with $h$ vertices and two tournaments
$T$ and $T^{\prime}$ defined on the same vertex set $V$ with $\mid
V\mid=n$. If $T$ and $T^{\prime }$ are $\{-p\}$-hypomorphic then for
each subset $X$ of $V$ of at most $p$ elements,
$n(T^{\prime},H;X)=n(T,H;X)$.
\end{corollary}

\begin{corollary}
\label{corollaire 1 pouz} \textrm{\cite{Pouz prop comb}} Consider
two tournaments $T=(V,A)$ and $T^{\prime}=(V,A^{\prime})$ and an
integer $p$ such that $0<p<\mid V\mid$. If $T$ and $T^{\prime}$ are
$\{p\}$-hypomorphic, then $T$ and $T^{\prime}$ are
$\{q\}$-hypomorphic for each $q\in\{1,\ldots ,min(p,\mid
V\mid-p)\}$.\\In particular, if $\mid V\mid \geq 6$ and $T$ and $T'$ are $\{-3\}$-hypomorphic,
then $T$ and $T'$ are $(\leq 3)$-hypomorphic.
\end{corollary}

Now, recall the following theorem, called the "Inversion Theorem",
which was obtained by A. Boussa\"{\i}ri, P. Ille, G. Lopez and S.
Thomass\'{e} (\cite{BILT1},\cite{BILT}).

\begin{theorem}
\label{theo inversion} (\textrm{\cite{BILT1},\cite{BILT}}) Given an
indecomposable tournament $T$ with at least $3$ vertices, the only
tournaments which are $\{3\}$-hypomorphic to $T$ are $T$ and
$T^{\ast}$.
\end{theorem}

The following corollary is a consequence of Theorem \ref{theo
inversion}.
\begin{corollary}
\label{rq quotient egaux ou invers} \textrm{\cite{BILT}} Let $T$ and
$T^{\prime}$ be two $\{3\}$-hypomorphic tournaments with at least
$3$ verices.

\begin{description}
\item[i)] $\mathcal{P}(T)=\mathcal{P}(T^{\prime})$.

\item[ii)] $T$ is strongly connected (resp. indecomposable) if and only if $T^{\prime}$ is strongly connected (resp. indecomposable).
\item[iii)] If $T$ is strongly connected, then the quotients
$T^{\prime }/\mathcal{P}(T)$ and $T/\mathcal{P}(T)$ are either equal
or dual.
\end{description}
\end{corollary}

From this corollary, we obtain the following remark.
\begin{remark}
\label{remark I interval of T donc I interval of T prime} Let $T$
and $T^{\prime}$ be two $\{3\}$-hypomorphic tournaments on a set $V$
with $\mid V\mid \geq 3$, and $I$ be a subset of $V$ such that
$T\left[  I\right]  $ is strongly connected. If $I$ is an interval
of $T$, then $I$ is an interval of $T'$.
\end{remark}

\bigskip

Given a tournament $T$ on a set $V$ and a subset $I$ of $V$, we
denote by $I^{+}_{T}$ (resp. $I^{-}_{T}$ ) the set of vertices $x\in
V\setminus I$ such that $I\longrightarrow x$ (resp.
$x\longrightarrow I$).

We complete this section by the following result.

\begin{proposition}
\label{propdegre}Let $T$ and $T^{\prime}$ be two tournaments on a
set $V$ \ with $\left\vert V\right\vert \geq 6$, and $I$ be an
interval of $T$ such that $\mid I\mid \geq 3$ and $T\left[  I\right]
$ is indecomposable.

\begin{description}
\item[i)] If $T$ and $T^{\prime}$ are $\{3,-2\}$-hypomorphic (resp. $\{-3\}$-hypomorphic) and $\left\vert
V\setminus I\right\vert \geq2$ (resp. $\left\vert V\setminus
I\right\vert \geq3$), then $T\left[  I\right]  \backsim
T^{\prime}\left[ I\right] $.

\item[ii)] If $T$ and $T^{\prime}$ are $\{3,-2\}$-hypomorphic (resp. $\{-3\}$-hypomorphic) and $\left\vert
V\setminus I\right\vert \geq3$ (resp. $\left\vert V\setminus
I\right\vert \geq4$), then $\left\vert I^{+}_{T}\right\vert
=\left\vert I^{+}_{T^{\prime}}\right\vert $ and $\left\vert
I^{-}_{T}\right\vert =\left\vert I^{-}_{T^{\prime}}\right\vert $.

\item[iii)] If $T$ and $T^{\prime}$ are $\{-3\}$-hypomorphic and $\left\vert
V\setminus I\right\vert \geq4$, then $T\left[  I\right]  \backsim
T^{\prime}\left[ I\right]  $, $I^{+}_{T}=I^{+}_{T^{\prime}}$ and
$I^{-}_{T}=I^{-} _{T^{\prime}}$.
\end{description}
\end{proposition}

\bpr

First, note that if $T$ and $T'$ are $\{-3\}$-hypomorphic, then $T$
and $T'$ are $\{3\}$-hypomorphic (by Corollary \ref{corollaire 1
pouz}). Moreover, as $T[I]$ is indecomposable (in particular, it is
strongly connected because $\mid I\mid \geq 3$) and $T[I]$ and
$T'[I]$ are $\{3\}$-hypomorphic, then by Corollary \ref{rq quotient
egaux ou invers}, $T'[I]$ is indecomposable and by Remark
\ref{remark I interval of T donc I interval of T prime}, $I$ is an
interval of $T'$.
\begin{description}
\item[i)] Let $a \neq b \in I$ and $J$ be a subset of $V$ containing $\{a,b\}$ such that $T[J]$ is indecomposable and $\mid I \mid = \mid J \mid$ and denote by $H$ the subtournament $T[I]$.
As $I$ is an interval of $T$, then $I \cap J$ is an interval of $T[J]$. However, $\{a,b\} \subset I \cap J$ and $T[J]$ is indecomposable, then $I \cap J = J$ and hence $I = J$ (because $\mid I\mid = \mid J\mid$). Thus, $I$ is the only subset $J$ of $V$ containing $\{a,b\}$ such that $T[J]$ is indecomposable and $\mid J\mid = \mid I\mid$. In particular $S(T,H;\{a,b\}) = \{I\}$ and then $n(T,H;\{a,b\}) = 1$.\\
By interchanging $T$ and $T'$ in the previous result, $I$ is the
only subset $J$ of $V$ containing $\{a,b\}$ such that $T'[J]$ is
indecomposable and $\mid I\mid =$ $\mid J\mid$. In particular,
$S(T',H;\{a,b\}) \subset \{I\}$. Lastly, as $T$ and $T'$ are
$\{-2\}$-hypomorphic (resp. $\{-3\}$-hypomorphic), $\mid V\mid \geq
6$ and $\mid I\mid \leq$

$\mid V\mid - 2$ (resp. $\mid I \mid \leq \mid V\mid - 3$), then by
Corollary \ref{corollaire 2 pouz}, $n(T,H;\{a,b\}) =
n(T',H;\{a,b\})$. As $n(T,H;\{a,b\}) = 1$, then $S(T',H;\{a,b\})
\neq \emptyset$ and so $S(T',H;\{a,b\}) = \{I\}$. Consequently,
$T'[I]\sim T[I]$.

\item[ii)] Let $a\neq b\in I$ and denote by $H$ a tournament with vertex set $I \cup \{u\}$
(where  $u \notin I $) such that $H[I] = T[I]$ and $I\longrightarrow
u$. Clearly, if $I^{+}_{T}\neq \emptyset$, then for each $x\in
I^{+}_{T}$, $I \cup \{x\} \in S(T,H;\{a,b\})$. Conversely, assume
that $S(T,H;\{a,b\}) \neq \emptyset$ and consider an element $J$ of
$S(T,H;\{a,b\})$. Let $f$ be an isomorphism from $H$ to $T[J]$ and
let $\alpha = f(u)$. As $I$ is the unique non trivial interval of
$H$ and $f(I) = J \setminus \{\alpha\}$, then $J \setminus
\{\alpha\}$ is the unique non trivial interval of $T[J]$. However,
$I \cap J$ is an interval of $T[J]$ and $\{a,b\} \subset I \cap J$,
then $I \cap J = J \setminus \{\alpha\}$ and hence $J \setminus
\{\alpha\} \subset I$. So, $J \setminus \{\alpha\} = I$ (because
$\mid I\mid = \mid J \setminus \{\alpha\} \mid$).
Thus, $J = I \cup \{\alpha\}$ and $\alpha \in I^{+}_{T}$.\\
We conclude that $S(T,H;\{a,b\}) = \{I \cup \{x\};\, x \in
I^{+}_{T}\}$ and hence,\\$n(T,H;\{a,b\}) =\mid I^{+}_{T}\mid$.\\
As $T'[I]\sim T[I]$ (by i)), then by interchanging $T$ and $T'$ in
the previous result, we deduce that $n(T',H;\{a,b\}) =$ $\mid
I^{+}_{T'} \mid$. Lastly, as $T$ and $T'$ are $\{-2\}$-hypomorphic
(resp. $\{-3\}$-hypomorphic) and $\mid I \cup \{u\} \mid \leq \mid V
\mid - 2$ (resp. $\mid I \cup \{u\}\mid \leq \mid V \mid - 3$), then
by Corollary \ref{corollaire 2 pouz}, $n(T,H;\{a,b\}) =
n(T',H;\{a,b\})$. Therefore, $\mid I^{+}_{T}\mid =\mid
I^{+}_{T'}\mid$ and hence, $\mid I^{-}_{T}\mid = \mid
I^{-}_{T'}\mid$.

\item[iii)] By ii), we have $\left\vert I^{+}_{T}\right\vert
=\left\vert I^{+}_{T^{\prime}}\right\vert $. Assume now that
$I^{+}_{T}\setminus I^{+}_{T^{\prime}}\neq\emptyset$ and let $x\in
I^{+}_{T}\setminus I^{+}_{T^{\prime}}$. As $I$ is an interval of
$T^{\prime}$, then $x\in I^{-}_{T^{\prime}}$. The tournaments $T-x$
and $T^{\prime}-x$ are $\{3,-2\}$-hypomorphic, $\mid V \setminus
\{x\}\mid \geq 6$, and $\left\vert (V\setminus \left\{ x\right\}
)\setminus I\right\vert \geq 3$, then by ii), $\left\vert
I^{+}_{T-x}\right\vert =\left\vert I^{+}_{T'-x}\right\vert $.
However, $\left\vert I^{+}_{T-x}\right\vert =\left\vert
I^{+}_{T}\right\vert -1\text{ and }\left\vert
I^{+}_{T'-x}\right\vert = \left\vert I^{+}_{T'}\right\vert $;
contradiction. It follows that $I^{+}_{T}\subset I^{+}_{T'}$ and
then $I^{+}_{T}=I^{+}_{T'}$. By duality, we obtain
$I^{-}_{T}=I^{-}_{T'}$. \end{description} \epr

\section{Proof of Theorem \ref{theorem principal moins 3 autodual}}

For the tournaments without diamonds, H. Bouchaala and Y. Boudabbous
\cite{Bouch et Boud} established the following result.

\begin{proposition}
\label{premiere prop de moin3 autod} \textrm{\cite{Bouch et Boud}}
Given a tournament $T$ without diamonds and which has at least $9$
vertices, $T$ is $\{-3\}$-self dual if and only if it is strongly
self dual.
\end{proposition}

We present now some results concerning tournaments embedding a
diamond.

\begin{remark}\label{remark delta I trois cycle}
A diamond $\delta$ has a unique non trivial interval $I$. Moreover, $\delta[I]$ is a $3$-cycle.
\end{remark}

\begin{lemma}\textrm{\cite{Bouch}}
\label{every vertex in a diamond} If $T=(V,A)$ is a tournament
embedding a diamond, then each vertex of $T$ is contained in at
least one diamond of $T$.
\end{lemma}

The following proposition was obtained by M. Sghiar in 2004. This
result plays an important role in the proof of Theorem \ref{theorem
principal moins 3 autodual}.

\begin{proposition} \textrm{\cite{Sghiar}}
\label{s il a inter a 2 donc non moins 3} Let $T$ be a tournament, with at
least $8$ vertices, embedding a diamond. If $T$ has an interval of cardinality
$2$, then $T$ is not $\{-3\}$-self dual.
\end{proposition}

For the proof of this proposition, we need some definitions and
notations. Given a tournament $T=(V,A)$, if $X=\{a,b,c,d\}$ is a
subset of $V$ such that $T[X]$ is a diamond and $T[\{a,b,c\}]$ is a
$3$-cycle, we say that $X$ is a diamond of $T$ of \emph{center} $d$
and cycle $\{a,b,c\}$. Let $x\neq y\in V$, we denote:

\begin{itemize}
\item $\delta_{T,\{x,y\}}^{+}$ (resp. $\delta_{T,\{x,y\}}^{-}$), the number of
positive (resp. negative) diamonds of $T$ whose cycle contains
$\{x,y\}$.

\item $C_{T,\{x,y\}}$, the set of elements $w$ of $V$ such that $T[\{x,y,w\}]$
is a $3$-cycle.

\item $\delta_{T,\{x,y,w\}}^{+}$ (resp. $\delta_{T,\{x,y,w\}}^{-}$), the
number of positive (resp. negative) diamonds of $T$ whose cycle is
$\{x,y,w\}$, where $w$ is an element of $C_{T,\{x,y\}}$.

\item $D_{T,\{x\}}^{+}(y)$ (resp. $D_{T,\{x\}}^{-}(y)$), the number of
positive (resp. negative) diamonds of $T$ passing by $x$ and whose center is $y$.

\item $D_{T,\{x,y\}}^{+}$ (resp. $D_{T,\{x,y\}}^{-}$), the number of positive
(resp. negative) diamonds of $T$ passing by $x$ and $y$.

\item $\delta_{T}^{+}(x)$ (resp. $\delta_{T}^{-}(x)$), the number of positive
(resp. negative) diamonds of $T$ whose center is $x$.
\end{itemize}

\begin{lemma}
\label{moins 3 donc every vertex is center}\textrm{\cite{Sghiar}}
Let $T = (V,A)$ be a $\{-3\}$-self dual tournament with at least $7$
vertices. If $T$ embeds a diamond, then every vertex of $T$ is the
center of at least one diamond of $T$.
\end{lemma}

\bpr Suppose for a contradiction that there exists a vertex $x$ of $T$ such that $\delta^{+}_{T}(x) =\delta^{-}_{T}(x) = 0$.
From Lemma \ref{every vertex in a diamond}, there exists a diamond $\sigma$ of $T$ containing $x$. By interchanging $T$ and $T^{\ast}$,
we can assume that $\sigma$ is a negative diamond. Let $y$ be the center of $\sigma$. So, $D^{-}_{T,\{x\}}(y) \neq 0$ and $D^{+}_{T,\{y\}}(x) = 0$.
If $C_{T,\{x,y\}} = \emptyset$, then $\delta^{+}_{T,\{x,y\}} = \delta^{-}_{T,\{x,y\}} = 0$. If $C_{T,\{x,y\}} \neq \emptyset$, then pick $w \in C_{T,\{x,y\}}$ and let $X = \{x,y,w\}$.
As $T$ and $T^{\ast}$ are $\{-3\}$-hypomorphic and $\mid V\mid \geq 7$, then from Corollary \ref{corollaire 2 pouz}, $n(T,\delta^{+};X) = n(T^{\ast},\delta^{+};X)$. So, $n(T,\delta^{+};X) = n(T,\delta^{-};X)$ and hence $\delta^{+}_{T,X} = \delta^{-}_{T,X}$. However, $\delta^{+}_{T,\{x,y\}} = \sum\limits_{w\in C_{T,\{x,y\}}}\delta_{T,\{x,y,w\}}
^{+}$ and $\delta_{T,\{x,y\}}^{-}=\sum\limits_{w\in C_{T,\{x,y\}}}
\delta_{T,\{x,y,w\}}^{-}$. Thus, $\delta^{+}_{T,\{x,y\}} = \delta^{-}_{T,\{x,y\}}$. On the other hand, we have $D^{-}_{T,\{x,y\}} = n(T,\delta^{-};\{x,y\}) = n(T^{\ast},\delta^{+};\{x,y\})$, $D^{+}_{T,\{x,y\}} = n(T,\delta^{+};\{x,y\})$ and from Corollary \ref{corollaire 2 pouz}, $n(T,\delta^{+};\{x,y\}) = n(T^{\ast},\delta^{+};\{x,y\})$, hence $D^{+}_{T,\{x,y\}} = D^{-}_{T,\{x,y\}}$. However, $D^{+}_{T,\{x,y\}} = D^{+}_{T,\{y\}}(x) + \delta^{+}_{T,\{x,y\}}$, $D^{-}_{T,\{x,y\}} = D^{-}_{T,\{x\}}(y) + \delta^{-}_{T,\{x,y\}}$ and $\delta^{+}_{T,\{x,y\}} = \delta^{-}_{T,\{x,y\}}$, thus, $D^{+}_{T,\{y\}}(x) = D^{-}_{T,\{x\}}(y)$; which contradicts the fact that $D^{-}_{T,\{x\}}(y) \neq 0$ and $D^{+}_{T,\{y\}}(x) = 0$.
\epr

\begin{lemma}
\label{moins 2 ou moins 3 inter 2}\textrm{\cite{Sghiar}} Consider a
$\{-2\}$-self dual (resp. $\{-3\}$-self dual) tournament $T = (V,A)$ with at
least $7$ (resp. $8$) vertices and two distinct vertices $a, b$ of $T$. If
$\{a,b\}$ is an interval of $T$ then $\delta^{+}_{T}(a) = \delta^{-}_{T}(a)$.
\end{lemma}

\bpr Let $H$ be the tournament obtained from one positive diamond by
dilating its center by a tournament with $2$ vertices. Let
$\Delta^{+}_{T}(a) = \{X\subset V;\, T[X]$ is a positive diamond
with center $a\}$. Let $X$ be an element of $\Delta^{+}_{T}(a)$. As
$\{a,b\}$ is an interval of $T$, so $\{a,b\} \cap X$ is an interval
of $T[X]$. Then, by Remark \ref{remark delta I trois cycle}, $b
\not\in X$ and $\{a,b\}$ is an interval of $T[X\cup \{b\}]$. Hence,
$X\cup \{b\} \in S(T,H;\{a,b\})$. Let's consider the map $f\, : \,
\Delta^{+}_{T}(a)\longrightarrow S(T,H;\{a,b\})$ defined by: for
each $X \in \Delta^{+}_{T}(a)$, $f(X) = X \cup \{b\}$. Clearly, $f$
is bijective and so $\delta^{+}_{T}(a) = n(T,H;\{a,b\})$. By
interchanging $T$ and $T^{\ast}$, we deduce that $\delta^{-}_{T}(a)
= \delta^{+}_{T^{\ast}}(a) = n(T^{\ast},H;\{a,b\})$. On the other
hand, as $T$ and $T^{\ast}$ are $\{-2\}$-hypomorphic (resp.
$\{-3\}$-hypomorphic) and $\mid V\mid \geq 7$ (resp. $\mid V\mid
\geq 8$), then from Corollary \ref{corollaire 2 pouz},
$n(T,H;\{a,b\}) = n(T^{\ast},H;\{a,b\})$. Thus, $\delta^{+}_{T}(a) =
\delta^{-}_{T}(a)$. \epr\\ 

\textbf{Proof of Proposition  $24$.} Assume by contradiction, that
$T$ is $\{-3\}$-self dual and has an interval $\{a,b\}$ with $a\neq
b$. By Lemma \ref{moins 3 donc every vertex is center}, $T$ has a
diamond $T\left[  X\right] $ with center $a$. Assume for example
that $T\left[  X\right]  $ is a positive diamond. Clearly, $b\notin
X$. Consider a vertex $x\in X\setminus \{a\}$. As $T-x$ (resp. $T$)
is $\{-2\}$-self dual (resp. $\{-3\}$-self dual) and $\{a,b\}$ is an
interval of $T-x$ (resp. $T$),
then, by Lemma \ref{moins 2 ou moins 3 inter 2}, $\delta_{T-x}^{+}%
(a)=\delta_{T-x}^{-}(a)$ (resp. $\delta_{T}^{+}(a)=\delta_{T}^{-}(a)$). So,
$0=\delta_{T}^{+}(a)-\delta_{T}^{-}(a)=\delta_{T-x}^{+}(a)+D_{T,\{x\}}%
^{+}(a)-\delta_{T-x}^{-}(a)=D_{T,\{x\}}^{+}(a)$; contradiction.\epr\\

Theorem \ref{theorem principal moins 3 autodual} is an
immediate consequence of Proposition  \ref{premiere prop de moin3
autod} and the below proposition.

\begin{proposition}
\label{deuxieme prop de moins 3 auto} Every decomposable tournament with at
least $8$ vertices embedding a diamond is not $\{-3\}$-self dual.
\end{proposition}

\bpr Let $T=(V,A)$ be a decomposable tournament with at least $8$
vertices embedding a diamond. Assume by contradiction that $T$ is
$\{-3\}$-self dual. By Proposition \ref{s il a inter a 2 donc non
moins 3}, $T$ has no interval of size $2$. Let $I$ be a minimal
(under the inclusion) non trivial interval of $T$. Clearly, the
subtournament $T\left[ I\right]  $ is indecomposable. And then it is
strongly connected. Firstly, assume that $\left\vert V\setminus
I\right\vert \geq4$. Let $z\in V\setminus I$ and suppose, for
example that $z\longrightarrow I$ in $T$. As $T$ and $T^{\ast}$ are
$\{-3\}$-hypomorphic, then from Proposition \ref{propdegre},
$z\longrightarrow I$ in $T^{\ast}$; contradiction. Secondly, assume
that $\left\vert V\setminus I\right\vert \leq3$. Let $k = \mid
V\setminus I\mid$. We have so $k \in \{1,2,3\}$ and $\mid I\mid \geq
8 - k$. As the subtournament $T[I]$ is strongly connected and $8 - k
\in \{5,6,7\}$, then from Lemma \ref{Moon}, there exists a subset
$X$ of $I$ such that $\mid X\mid = 4 - k$ and $T[I] - X$ is strongly
connected. Let $Y$ be a subset of $V\setminus I$ such that $\mid
Y\mid = k - 1$. Clearly, the subtournament $T - (X \cup Y)$ is not
self dual; contradicts the fact that $\mid X \cup Y\mid = 3$. \epr

\section{Proof of Theorem \ref{Theo moins 3 reconstruction}}

The proof of Theorem \ref{Theo moins 3 reconstruction}  is based on the next result.

\begin{proposition}
\label{prop2 de moin3 recons} Let $T$ be a strongly connected and
decomposable tournament on a set $V$ with $\mid V\mid = n\geq 9$,
which is not almost transitive. If $T^{\prime}$ is a tournament
$\{-3\}$-hypomorphic to $T$, then the following assertions hold.

\begin{enumerate}
\item $\mathcal{P}(T^{\prime}) = \mathcal{P}(T)$ and $T^{\prime}%
/\mathcal{P}(T) = T/\mathcal{P}(T)$.

\item If there exists $X\in\mathcal{P}(T)$ such that $T^{\prime}%
[X]\not \sim T[X]$, then $\mid\mathcal{P}(T)\mid=3$ and $\mid X\mid=n-2$.

\item If for each $X \in \mathcal{P}(T)$, $\mid X\mid \leq n - 3$, then for each $X \in \mathcal{P}(T)$, $T'[X]\sim T[X]$ and in particular $T' \sim T$.
\end{enumerate}
\end{proposition}

For the proof of the last proposition, we use the following remark and lemma.

\begin{remark}\label{remark isomorphi hereditaire}
Let $T$ and $T'$ be two tournaments on a set $V$ with $\mid V\mid \geq 3$ and $\Gamma$ be a common interval partition of $T$ and $T'$ such that $T/\Gamma = T'/\Gamma$. Given a non empty subset $A$ of $V$ and let $\Gamma_{A} = \{X \cap A;\, X\in \Gamma$ and $ X\cap A \neq \emptyset\}$. Then $\Gamma_{A}$ is a common interval partition of $T[A]$ and $T'[A]$ and $T[A]/\Gamma_{A} = T'[A]/\Gamma_{A}$. Suppose that for each $Y \in \Gamma_{A}$, there exists an isomorphism $\varphi_{Y}$ from $T[Y]$ onto $T'[Y]$ and consider the map $f\, :\, A\longrightarrow A$ defined by: for each $x \in A$, $f(x) = \varphi_{Y}(x)$ where $Y$ is the unique element of $\Gamma_{A}$ such that $x \in Y$. Then, $f$ is an isomorphism from $T[A]$ onto $T'[A]$. In particular, if for each $X\in \Gamma$, $T[X]$ and $T'[X]$ are hereditarily isomorphic, then $T$ and $T'$ are hereditarily isomorphic.
\end{remark}

\begin{lemma}
\label{bouss gene}\textrm{\cite{Boua boud}} Let $T=(V,A)$ and $T^{\prime
}=(V',A^{\prime})$ be two isomorphic tournaments, $f$ be an isomorphism from
$T$ onto $T^{\prime}$, $i\in V$ and $R_{i}$ (resp. $R_{i}^{\prime}$) be a
tournament defined on a vertex set $I_{i}$ (resp. $I_{i}^{\prime}$) disjoint
from $V$ (resp. $V^{\prime}$). Let $R$ (resp. $R^{\prime}$) be the tournament
obtained from $T$ (resp. $T^{\prime}$) by dilating the vertex $i$ (resp.
$f(i)$) by $R_{i}$ (resp. $R_{i}^{\prime}$). Then $R\sim R^{\prime}$ if and
only if $R_{i}\sim R_{i}^{\prime}$.
\end{lemma}

Note that this lemma is a simple generalization of a result
communicated by A. Boussa\"{\i}ri, and on which $V' = V$ and $f=id_{V}$.\\

\newpage
\textbf{Proof of Proposition \ref{prop2 de moin3 recons}.}
\begin{enumerate}
\item By Corollary \ref{corollaire 1 pouz}, $T$ and $T^{\prime}$ are $(\leq
3)$-hypomorphic. So, from Corollary \ref{rq quotient egaux ou
invers},
$\mathcal{P}(T)=\mathcal{P}(T^{\prime})$ and $T^{\prime}/\mathcal{P}%
(T)=T/\mathcal{P}(T)$ or $T^{\prime}/\mathcal{P}(T)=T^{\ast}/\mathcal{P}%
(T)$.\newline Assume by contradiction that $T^{\prime}/\mathcal{P}(T)=T^{\ast
}/\mathcal{P}(T)$. In this case, we are going to show that for every
$X\in\mathcal{P}(T)$, $T[X]$ is transitive, and thus by Remark \ref{remark isomorphi hereditaire}, $T^{\prime}$ is
hereditarily isomorphic to $T^{\ast}$. Since $T^{\prime}$ is $\{-3\}$%
-hypomorphic to $T$, then $T$ is $\{-3\}$-self dual. By Theorem
\ref{theorem principal moins 3 autodual}, the tournament $T$ is
almost transitive; which contradicts the hypothesis. For that,
proceed by contradiction and consider an element $X$ of
$\mathcal{P}(T)$ and a subset $\{\alpha,\beta,\gamma\}$ of $X$ such
that $T[\{\alpha,\beta,\gamma\}]$ is a $3$-cycle. As $T$ is strongly
connected, there is $a\in V\setminus X$ such that $X\longrightarrow
a$. From Corollary \ref{corollaire 2 pouz}, $n(T,\delta
^{+};\{\alpha,\beta,a\})=n(T^{\prime},\delta^{+};\{\alpha,\beta,a\})$.
Moreover, $n(T,\delta^{+};\{\alpha,\beta,a\})\neq0$, because
$T[\{\alpha ,\beta,\gamma,a\}]\sim\delta^{+}$, then there exists a
subset $K$ of $V$ such that $\{\alpha,\beta,a\}\subset K$ and
$T^{\prime}[K]$ $\sim\delta^{+}$. Hence, $\mid K\cap X\mid=3$,
because otherwise $K\cap X$ is an interval with two elements of
$T^{\prime}[K]$; which contradicts Remark \ref{remark delta I trois
cycle}. So, $T^{\prime}[K]$ is written: $(K\cap X)\longleftarrow a$;
which contradicts the fact that $T^{\prime}[K]\sim\delta^{+}$.
\item By 1, we have $\mathcal{P}(T') = \mathcal{P}(T)$ and $T^{\prime}/\mathcal{P}(T)=T/\mathcal{P}(T)$. We
distinguish the following two cases.\begin{itemize}
\item If for every $X\in\mathcal{P}(T)$, $\mid X\mid\leq n-\left\vert
\mathcal{P}(T)\right\vert -2$.\newline Let $X\in\mathcal{P}(T)$ and
let $H$ be the tournament obtained from $T/\mathcal{P}(T)$ by
dilating the vertex $X$ by $T[X]$. Assume that $\mid X\mid \geq 2$
and consider a subset $A$ of $X$ with $2$ elements. Consider a
subset $B$ of $V$ containing $X$ such that: $\forall\, Y\in
\mathcal{P}(T)\setminus \{X\}$, $\mid B\cap Y\mid = 1$. Clearly,
$B\in S(T,H;A)$ and hence $n(T,H;A) \neq \emptyset$. From Corollary
\ref{corollaire 2 pouz}, $n(T^{\prime},H;A)=n(T,H;A)$, then
$n(T^{\prime},H;A)\neq0$. So, there exists a subset $K$ of $V$ such
that: $A\subset K$ and $T^{\prime}[K]\sim H$. Then
$\mathcal{P}(T'[K])$ have a unique element $J$ non reduced to a
singleton.
Clearly, $T'[J]\sim T[X]$, in particular, $\mid J\mid = \mid X\mid$.\\
Let $P_{1} = \{Y\in \mathcal{P}(T);\, \mid Y\cap K\mid \geq 2\}$ and
$P_{2} = \{Y\in \mathcal{P}(T);\\ \mid Y\cap K\mid = 1\}$. The set
$K$ is the union of the two disjoint sets $K_{1} =
\bigcup\limits_{Y\in P_{1}} K \cap Y$ and $K_{2} =
\bigcup\limits_{Y\in P_{2}} K \cap Y$. As $A \subset X\cap K$, then
$X \in P_{1}$ and so, $P_{2}\subset (\mathcal{P}(T) \setminus
\{X\})$, in particular, $\mid P_{2}\mid \leq \mid \mathcal{P}(T)\mid
- 1$. For all $Y\in P_{1}$, $Y \cap K$ is a non trivial interval of
$T'[K]$, then $Y\cap K \subset J$, so $K_{1}\subset J$ and thus
$(K\setminus J) \subset K_{2}$.\\ Thus, $\mid \mathcal{P}(T)\mid - 1
= \mid K \setminus J\mid \leq \sum\limits_{Y\in P_{2}}\mid K\cap
Y\mid = \mid P_{2}\mid\leq \mid \mathcal{P}(T)\mid - 1$. So, $\mid
P_{2}\mid = \mid \mathcal{P}(T)\mid - 1$ and then $P_{2} =
\mathcal{P}(T)\setminus \{X\}$ and $P_{1} = \{X\}$. Thus, $K =
(X\cap K) \cup K_{2}$. So, $\mid X\cap K\mid = \mid K\mid - \mid
K_{2}\mid =$ $\mid K\mid - \mid P_{2}\mid = \mid K\mid - \mid
\mathcal{P}(T)\mid + 1 = \mid J\mid$. As in addition $X\cap K\subset
J$, then, $X \cap K = J$. So, $J\subset X$ and then $J = X$ because
$\mid J\mid = \mid X\mid$. Consequently, $T'[X]\sim T[X]$.
\item If there exists an element $X$ of $\mathcal{P}(T)$ such that $\mid
X\mid>n-\left\vert \mathcal{P}(T)\right\vert -2$.\newline  In this
case, it is clear that for every $Y\in\mathcal{P}(T)\setminus
\{X\}$, $\mid Y\mid\leq3$, so $T^{\prime}[Y]$ and $T[Y]$ are
hereditarily isomorphic. Suppose that $\mid V\setminus X\mid\geq3$.
Let $x$ be an element of $X$ and $B$ be a subset of $V\setminus X$
such that $\mid B\mid=3$. Denote by $V_{(x,B)}$ the set
$(V\setminus(X\cup B))\cup\{x\}$ and by $T_{(x,B)}$ (resp.
$T_{(x,B)}^{\prime}$) the subtournament $T[V_{(x,B)}]$ (resp.
$T^{\prime}[V_{(x,B)}]$). The tournament $T-B$ (resp.
$T^{\prime}-B$) is obtained from the tournament $T_{(x,B)}$ (resp.
$T_{(x,B)}^{\prime}$) by dilating the vertex $x$ by $T[X]$ (resp.
$T^{\prime}[X]$). Moreover, by Remark \ref{remark isomorphi
hereditaire}, there exists an isomorphism $g$ from $T_{(x,B)}$ onto
$T_{(x,B)}^{\prime}$ such that $g(x)=x$. As in addition, $T-B$ and
$T^{\prime}-B$ are isomorphic, then, from Lemma \ref{bouss gene},
$T^{\prime }[X]$ is isomorphic to $T[X]$; which permits to conclude.
\end{itemize}
\item Is a direct consequence of 2.
\end{enumerate}
\hskip 12cm $\square$

\begin{lemma}
\label{fc dec 8 vertx est moins 2 3 reco} Every strongly connected
and decomposable tournament, which has $8$ vertices, is
$\{-2,-3\}$-reconstructible.
\end{lemma}

\bpr Let $H$ be a strongly connected and decomposable tournament
defined on a vertex set $X$ with $\mid X\mid=8$ and let $H^{\prime}$
be a tournament $\{-2,-3\}$-hypomorphic to $H$. The tournaments $H$
and $H^{\prime}$ are $\{3\}$-hypomorphic, by Corollary
\ref{corollaire 1 pouz}, so they are $\{3,5,6\}$-hypomorphic.
Besides, by Corollary \ref{rq quotient egaux ou invers},
$\mathcal{P}(H) = \mathcal{P}(H')$, and $H'/\mathcal{P}(H) =
H/\mathcal{P}(H)$ or $H'/\mathcal{P}(H) = H^{\ast}/\mathcal{P}(H)$.
Let's put $Q=\mathcal{P}(H)$, and discuss according to its cardinal.

\begin{itemize}
\item If $\mid Q\mid>3$.\newline If for every $Z\in Q$, $\mid Z\mid\leq3$. By
$(\leq3)$-hypomorphy, $H^{\prime}[Z]\sim H[Z];\,\forall Z\in Q$. If
$H^{\prime}/Q=H/Q$, clearly $H^{\prime}\sim H$. Suppose hence that
$H^{\prime}/Q = H^{\ast}/Q$. In this case, $H^{\prime}$ is
hereditarily isomorphic to $H^{\ast}$. Then, $H$ is $\{-3\}$-self
dual. From Proposition \ref{deuxieme prop de moins 3 auto}, $H$ is
without diamonds. By Remark \ref{tournaments with 4 vertices}, iii),
$H$ is $\{4\}$-self dual. Thus, $H$ is $\{4,5,6\}$-self dual. So,
$H$ is $(\leq6)$-self dual and it is thus self dual, by the $(\leq
6)$-reconstruction of tournaments with at least $6$ vertices
\cite{Lopez indéf}. It follows that $H^{\prime}\sim H$. \newline If
there exists $Y\in Q$ such that $\mid Y\mid\geq4$. In this case, as
$\mid X\mid=8$, $\mid Q\mid>3$ and every tournament with $4$
vertices is decomposable, then $\mid Q\mid=5$, $\mid Y\mid=4$ and
for every $Z\in Q\setminus \{Y\}$, $\mid Z\mid=1$.

\begin{itemize}
\item If $H^{\prime}/Q=H/Q$. For $z\in X\setminus Y$, as $\mid Y\cup\{z\}\mid=\mid
X\mid-3$, then, $H[Y\cup\{z\}]\sim H^{\prime}[Y\cup\{z\}]$. It
follows that $H^{\prime}[Y]\sim H[Y]$, and then $H^{\prime}\sim H$.

\item If $H^{\prime}/Q = H^{\ast}/Q$.

\begin{itemize}
\item If $H[Y]$ is not a diamond. In this case, clearly $H^{\prime}[Y]$ is
hereditarily isomorphic to $H^{\ast}[Y]$ and hence $H'$ and $H^{\ast}$ are hereditarily isomorphic. As $H'$ is $\{-3\}$-hypomorphic to $H$, then $H$ is
$\{-3\}$-self dual. From Proposition \ref{deuxieme prop de moins 3
auto}, $H$ is then without diamonds, and it is clearly self dual and
thus, $H^{\prime}\sim H$.

\item If $H[Y]$ is a diamond. In this case, if $H^{\prime}[Y]$ is not
isomorphic to $H[Y]$, then $H^{\prime}$ is hereditarily isomorphic
to $H^{\ast}$ and so, $H$ is $\{-3\}$-self dual; which contradicts
Proposition \ref{deuxieme prop de moins 3 auto}. So,
$H^{\prime}[Y]\sim H[Y]$. For $z\in X\setminus Y$, it is easy to
verify that $H^{\prime}[Y\cup\{z\}]$ is not isomorphic to
$H[Y\cup\{z\}]$; which contradicts the $\{-3\}$-hypomorphy between
$H^{\prime }$ and $H$.
\end{itemize}
\end{itemize}

\item If $\mid Q\mid=3$. We distinguish the following two sub-cases.

\begin{itemize}
\item If there exists an unique $Y\in Q$ such that $\mid Y\mid>1$.\newline In
this case, $Q$ is written: $Q=\{\{a\},\{b\},Y\}$ where $\mid
Y\mid=6$ and $a\longrightarrow Y\longrightarrow b$ in $H$. As $\mid
Y\mid=6$, we have: $H^{\prime}[Y]\sim H[Y]$ by $\{-2\}$-hypomorphy.
Thus, $H^{\prime}\sim H$.

\item If there are $Y\neq Z\in Q$ such that $\min(\mid Y\mid,\mid Z\mid
)>1$.\newline In this case, we can write: $Q=\{Y_{1},Y_{2},Y_{3}\}$
with: $\mid Y_{3}\mid\leq\mid Y_{2}\mid\leq$ $\mid Y_{1}\mid$, $\mid
Y_{2}\mid\geq2$ and $\mid Y_{1}\mid\geq3$. By considering a subset
$A_{1}$ (resp. $A_{2}$) with $3$ (resp. $2$) elements of $Y_{1}$
(resp. $Y_{2}$) and an element
$y_{3}$ of $Y_{3}$, we see that the isomorphy between $H[A_{1}\cup A_{2}%
\cup\{y_{3}\}]$ and $H^{\prime}[A_{1}\cup A_{2}\cup\{y_{3}\}]$ requires that:
$H^{\prime}/Q=H/Q$. So, if $\mid Y_{1}\mid=3$, then for every $i\in\{1,2,3\}$,
$H^{\prime}[Y_{i}]\sim H[Y_{i}]$ and then $H^{\prime}\sim H$. Suppose thus
that $\mid Y_{1}\mid\geq4$. As $\mid X\mid=8$, then $\mid Y_{1}\mid\in\{4,5\}$
and $\mid Y_{2}\mid\leq3$. By considering an element $y_{2}$ of $Y_{2}$, we
see that the isomorphy between $H[Y_{1}\cup\{y_{2}\}]$ and $H^{\prime}%
[Y_{1}\cup\{y_{2}\}]$ requires the isomorphy between $H^{\prime}[Y_{1}]$ and
$H[Y_{1}]$. It follows that $H^{\prime}\sim H$.
\end{itemize}
\epr
\end{itemize}

\textbf{Proof of Theorem \ref{Theo moins 3 reconstruction}.}
Consider a decomposable tournament $T$ defined on a vertex set $V$
with $\mid V\mid=n\geq9$, a tournament $T^{\prime}$
$\{-3\}$-hypomorphic to $T$.
By Corollary \ref{corollaire 1 pouz}, $T'$ is $(\leq 3)$-hypomorphic to $T$. So, from Corollary \ref{rq quotient egaux ou invers}, $\mathcal{P}(T) = \mathcal{P}(T')$. In particular, if $T$ is strongly connected, then $\widetilde{P}(T) = \widetilde{P}(T')$.\\
If $T$ is almost transitive, then clearly $T' \sim T$. Let's suppose that $T$ is not almost transitive.
For the proof, we distinguish the following two cases.

\begin{description}
\item[ $\ast$] \textbf{ }$T$\textbf{ } is strongly connected.
\end{description}

In this case, from Proposition \ref{prop2 de moin3 recons}, $T^{\prime}/\mathcal{P}(T)=T/\mathcal{P}
(T)$ and we can suppose that $\mid\mathcal{P}(T)\mid=3$ and there exists
$X\in\mathcal{P}(T)$ such that $\mid X\mid=n-2$. Let $\mathcal{P}(T)=\{X,\{a\},\{b\}\}$ where
$X\longrightarrow a\longrightarrow b\longrightarrow X$ in $T$. We verify easily that $T'[X]$ and $T[X]$ are $\{-1,-2,-3\}$-hypomorphic (because $T'$ and $T$ are $\{-3\}$-hypomorphic).
We consider the following three cases.

\begin{itemize}
\item  $T[X]$ is non-strongly connected.

As the tournaments $T^{\prime}[X]$ and $T[X]$ are
$\{-1\}$-hypomorphic, then they are isomorphic, since the
non-strongly connected tournaments with at least $5$ vertices are
$\{-1\}$-reconstructible \cite{Har et Palm}.

\item $T[X]$ is strongly connected and decomposable.

Let $Q=\mathcal{P}(T[X])$. Assume that there exists $Y\in Q$ such
that $\mid Y\mid=$ \\ $\mid X\mid-2$. In this case, we have $\mid
Q\mid=3$ and $T[X]/Q$ is a $3$-cycle. Moreover, as $T^{\prime}[X]$
and $T[X]$ are $(\leq3)$-hypomorphic, then by Corollary \ref{rq
quotient egaux ou invers},
$\mathcal{P}(T^{\prime}[X])=\mathcal{P}(T[X])$ and $T^{\prime}[X]/\mathcal{P}%
(T[X])=T[X]/\mathcal{P}(T[X])$ or
$T^{\prime}[X]/\mathcal{P}(T[X])=T^{\ast }[X]/\mathcal{P}(T[X])$. As
in addition, $T'[Y] \sim T[Y]$ (because $T'[X]$ and $T[X]$ are
$\{-2\}$-hypomorphic), then $T^{\prime}[X]\sim T[X]$. Thus, clearly $T' \sim T$. \\
Now, suppose that for every $Y\in Q$, $\mid Y\mid<\mid X\mid -2$. In
this case, $T[X]$ is not almost transitive. If $\mid X\mid \geq9$,
as $T'[X]$ and $T[X]$ are $\{-3\}$-hypomorphic, then by Proposition
\ref{prop2 de moin3 recons}, $T'[X]\sim T[X]$ and hence, clearly
$T'\sim T$. If $\mid X\mid=7$. As $T^{\prime}[X]$ and $T[X]$ are
$\{-1,-2,-3\}$-hypomorphic, then they are $(\leq6)$-hypomorphic, and
thus $T^{\prime}[X]\sim T[X]$ (by \cite{Lopez indéf}). If $\mid
X\mid=8$. As $T^{\prime}[X]$ and $T[X]$ are $\{-2,-3\}$-hypomorphic,
then, by Lemma \ref{fc dec 8 vertx est moins 2 3 reco},
$T^{\prime}[X]\sim T[X]$.

\item If $T[X]$ is indecomposable.\newline In this case, from Theorem \ref{theo inversion},
$T'[X] = T[X]$ or $T'[X] = T^{\ast}[X]$ (because $T'[X]$ and $T[X]$
are $(\leq 3)$-hypomorphic). If $T'[X] = T[X]$, then clearly $T'\sim
T$. If $T'[X] = T^{\ast}[X]$, then $T[X]$ is $\{-1,-2,-3\}$-self
dual (because $T'[X]$ and $T[X]$ are $\{-1,-2,-3\}$-hypomorphic). We
obtain: If $T[X]$ is self dual, then $T'[X]\sim T[X]$ and clearly
$T' \sim T$. If $T[X]$ is not self dual, then $T \in
C_{3}(I_{n-2,\{-1,-2,-3\}})\subset \Omega_{n} $ and clearly $T'
\not\sim T$.
\end{itemize}

\begin{description}
\item[$\ast$] $T$ is non-strongly connected.

Observe that if $T$ is a transitive tournament, then $T'$ is also
transitive (by $(\leq 3)$-hypomorphy) and the result is obvious.
Suppose then that $T$ is a non-strongly connected tournament which
is not transitive. The result follows from the following five facts.
\item \textbf{Fact 1}.  Let $X$ be an element of $\widetilde{\mathcal{P}}(T)$
such that $T[X]$ is strongly connected with $\mid X\mid\geq3$ and
let $a$ be an element of $V\setminus X$. As $\mathcal{P}(T') =
\mathcal{P}(T)$, then $T'[X]$ is strongly connected, $X \in
\widetilde{P}(T')$ and we have: \newline$a\longrightarrow X$ in $T$
if and only if $a\longrightarrow X$ in $T^{\prime}$.

\item Indeed :\newline Let $\{\alpha,\beta,\gamma\}$ be a subset of $X$ such
that $T[\{\alpha,\beta,\gamma\}]$ is a $3$-cycle and suppose, for
example, that $X\longrightarrow a$ in $T$. From Corollary
\ref{corollaire 2 pouz}, $n(T,\delta^{+};\{\alpha,\beta,a\})$
$=n(T^{\prime},\delta^{+};\{\alpha,\beta,a\})$. As
$n(T,\delta^{+};\{\alpha,\beta,a\})\neq0$, thus there exists a
subset $K$ of $V$ such that $\{\alpha,\beta,a\}\subset K$ and
$T^{\prime}[K]\sim \delta^{+}$. Hence, $\mid K\cap X\mid=3$, because
otherwise $K\cap X$ is an interval with two elements of $T'[K]$;
which contradicts Remark \ref{remark delta I trois cycle}. As $a\in
V\setminus X$ and $K\setminus \{a\}\subset X$, then $K\setminus
\{a\}$ is an interval of $T^{\prime}[K]$ and thus $T^{\prime}[K]$ is
a diamond of center $a$. So, $X\longrightarrow a$ in $T^{\prime}$.

\item \textbf{Fact 2}.  $\widetilde{\mathcal{P}}(T^{\prime})=\widetilde
{\mathcal{P}}(T)$.

\item Indeed :\\
Consider an element $Y$ of $\widetilde{\mathcal{P}}(T)$. We
distinguish the following two cases.

\begin{itemize}
\item If $T[Y]$ is strongly connected and $\mid Y\mid\geq3$. In this case,
$Y\in\mathcal{P}(T)$. Then $Y\in\mathcal{P}(T^{\prime})$ and hence,
$Y\in\widetilde{\mathcal{P}}(T^{\prime})$.

\item If $T[Y]$ is a transitive. In this case, there exists an element $Z$ of
$\widetilde{\mathcal{P}}(T^{\prime})$ such that $Y\subset Z$,
because otherwise, $Y$ admits a partition $\{Y_{1},Y_{2}\}$ such
that there is $K\in\mathcal{P}(T^{\prime})$ with $\mid K\mid\geq3$,
$T^{\prime}[K]$ is strongly connected and in $T^{\prime}$ we have
$Y_{1}\longrightarrow K\longrightarrow Y_{2}$; which contradicts the
Fact 1. While exchanging the roles of $T$ and $T^{\prime}$, we can
hence deduce that $Y=Z$ and then, $Y\in \widetilde{\mathcal{P}}(T')$
.
\end{itemize}

\item \textbf{Fact 3. }$T^{\prime}/\widetilde{\mathcal{P}}(T)=T/\widetilde
{\mathcal{P}}(T)$.

\item Indeed:

Proceed by the absurd and suppose that there exist two distinct
elements $X$ and $Y$ of $\widetilde{\mathcal{P}}(T)$ such that
$X\longrightarrow Y$ in $T$ and $Y\longrightarrow X$ in
$T^{\prime}$. From the Fact 1,  $T[X]$ and $T[Y]$ are transitive.
So, $X$ and $Y$ are not consecutive in $T/\widetilde
{\mathcal{P}}(T)$. Thus, there exists an element $Z$ of $\widetilde
{\mathcal{P}}(T)$ such that $T[Z]$ is strongly connected, $\mid
Z\mid\geq3$ and $X\longrightarrow Z\longrightarrow Y$ in $T$. So, by
Fact 1, $X\longrightarrow Z\longrightarrow Y$ in $T^{\prime}$ and
then $X\longrightarrow Y$ in $T^{\prime}$; which is absurd.

\item \textbf{Fact 4}. If $T \not\in \Omega_{n}$, then for all $X\in\widetilde{\mathcal{P}}(T)$, $T^{\prime
}[X]\sim T[X]$.

\item Indeed:

Suppose that $T \not\in \Omega_{n}$ and consider an element $X$ of $\widetilde{\mathcal{P}}(T)$. As $T$ and
$T^{\prime}$ are $(\leq3)$-hypomorphic, we can assume that $\mid X\mid\geq4$ and $T[X]$ is strongly
connected.

We distinguish the following  cases.

\begin{itemize}
\item If $\mid X\mid\leq n-3$.\newline Consider $H=T[X]$ and $A\subset X$ such that $\mid A\mid=2$. From Corollary
\ref{corollaire 2 pouz}, $n(T,H;A)=n(T^{\prime},H;A)$. As
$n(T,H;A)\neq0$, then $n(T^{\prime},H;A)\neq0$. So, there exists a
subset $K$ of $V$ such that $A\subset K$ and $T^{\prime}[K]\sim H$.
We have then $K=X$, because otherwise, there exists
$Y\in\widetilde{\mathcal{P}}(T)\setminus \{X\}$ such that $K\cap
Y\neq\emptyset$ and hence $T^{\prime }[K]$ is non-strongly connected
(because $K\cap X$ is also non empty); which is absurd. So,
$T^{\prime}[X]\sim H=T[X]$.

\item If $\mid X\mid=n-1$.\newline Let $\{a\}=V\setminus X$ and suppose, for example,
that $X\longrightarrow a$ in $T$. As $\mid X\mid-2=n-3$, then
$T^{\prime}[X]$ and $T[X]$ are $\{-2\}$-hypomorphic. Now we shall
prove that these two tournaments are $\{-3\}$-hypomorphic. For that,
consider a subset $A$ of $X$ such that $\mid A\mid=3$. It is clear
that $T^{\prime}-A$ and $T-A$ are isomorphic and in these two
tournaments, we have $(X\setminus A)\longrightarrow a$. So,
$T^{\prime}[X\setminus A]\sim T[X\setminus A]$. Hence,
$T^{\prime}[X]$ and $T[X]$ are
$\{-3\}$-hypomorphic. So, $T^{\prime}[X]$ and $T[X]$ are $\{-2,-3\}$%
-hypomorphic.
\begin{itemize}
\item If $T[X]$ is decomposable. Let $Q=\mathcal{P}(T[X])$.
If there exists $ $ $ $ $Y\in Q$ such that $\mid Y\mid=\mid
X\mid-2$. In this case, $\mid Q \mid=3$ and $T[X]/Q$ is a $3$-cycle.
As $\mid Y\mid=n-3$, then $T^{\prime}[Y]\sim T[Y]$. As in addition,
$T^{\prime}[X]$ and $T[X]$ are $(\leq3)$-hypomorphic, then, by
Corollary \ref{rq quotient egaux ou invers},
  $\mathcal{P}(T^{\prime
}[X])=Q$ and $T^{\prime}[X]/Q$ $ =T[X]/ Q$ or $T^{\prime}[X]/ Q$
$=T^{\ast}[X]/ Q$. So we deduce immediately that $T^{\prime}[X]\sim
T[X]$. Now suppose that for every $Y\in Q$, $\mid Y\mid<\mid
X\mid-2$. In this case, $T[X]$ is not almost transitive. Distinguish
the following two cases. If $\mid X\mid\geq9$. In this case, by
applying Proposition \ref{prop2 de moin3 recons} for the tournaments
$T[X]$ and $T'[X]$, we obtain $T^{\prime}[X]\sim T[X]$. If $\mid
X\mid=8$. As $T^{\prime}[X]$ and $T[X]$ are $\{-2,-3\}$-hypomorphic,
then, by Lemma \ref{fc dec 8 vertx est moins 2 3 reco},
$T^{\prime}[X]\sim T[X]$.
\item If $T[X]$ is indecomposable. We have $T^{\prime}[X]$ and
$T[X]$ are $(\leq 3)$-hypomorphic. So, by Theorem \ref{theo
inversion}, $T^{\prime}[X]=T[X]$ or $T^{\prime}[X]=T^{\ast}[X]$.
Consider then the case where $T^{\prime }[X]=T^{\ast}[X]$. If $T[X]$
is self dual, then $T^{\prime}[X]\sim T[X]$. Suppose hence that
$T[X]$ is not self dual. As $T[X]$ is $\{-2,-3\}$-self dual
(because $T^{\prime}[X]$ and $T[X]$ are $\{-2,-3\}$%
-hypomorphic), then $T\in O_{2}(I_{n-1,\{-2,-3\}}) \subset
\Omega_{n}$; which is absurd.
\end{itemize}
\item If $\mid X\mid=n-2$.\newline Let $\{a,b\}=V\setminus X$ and suppose,
for example, that $X\longrightarrow a$ in $T$. As $\mid
X\mid-1=n-3$, then $T^{\prime}[X]$ and $T[X]$ are
$\{-1\}$-hypomorphic. Now, we shall prove that $T^{\prime}[X]$ and
$T[X]$ are $\{-2,-3\}$-hypomorphic. For that, for every
$i\in\{2,3\}$, let $A_{i}$ be a subset of $X$ such that $\mid
A_{i}\mid=i$. We have $T-(A_{2}\cup\{a\})\sim
T^{\prime}-(A_{2}\cup\{a\})$. From Fact 3, $(X\setminus
A_{2})\longrightarrow b$ in $T$ if and
only if $(X\setminus A_{2})\longrightarrow b$ in $T^{\prime}$. So, $T^{\prime}%
[X\setminus A_{2}]\sim T[X\setminus A_{2}]$ and thus,
$T^{\prime}[X]$ and $T[X]$ are $\{-2\}$-hypomorphic. Furthermore, as
$T^{\prime}-A_{3}\sim T-A_{3}$, then by
Fact 3, we can see that $T^{\prime}%
[X\setminus A_{3}]\sim T[X\setminus A_{3}]$. So, $T^{\prime}[X]$ and $T[X]$ are $\{-1,-2,-3\}$%
-hypomorphic.
\begin{itemize}
\item If $T[X]$ is decomposable. Pose
$Q=\mathcal{P}(T[X])$. If there exists $Y\in Q$ such that $\mid
Y\mid=\mid X\mid-2$. In this case, $\mid Q \mid=3$ and $T[X]/Q$ is a
$3$-cycle. Besides, by $\{-3\}$-hypomorphy, the two subtournaments
$T^{\prime}[Y\cup\{a\}]$ and $T[Y\cup\{a\}]$ are isomorphic (because
$\mid Y\cup\{a\}\mid=n-3$). As in addition, $Y\longrightarrow a$ in
both these tournaments, then $T^{\prime}[Y]\sim T[Y]$. As in
addition, $T^{\prime}[X]$ and $T[X]$ are $(\leq3)$-hypomorphic,
then, by Corollary \ref{rq quotient egaux ou invers}, we can see
that $T^{\prime}[X]\sim T[X]$. Now, suppose that for every $Y\in Q$,
$\mid Y\mid<\mid X\mid-2$. In this case, $T[X]$ is not almost
transitive. We distinguish the following three cases.\\ If $\mid
X\mid\geq9$, we conclude by Proposition \ref{prop2 de moin3
recons}.\\ If $\mid X\mid=8$, we conclude by Lemma \ref{fc dec 8
vertx est moins 2 3 reco}.\\ If $\mid X\mid=7$, as $T^{\prime}[X]$
and $T[X]$ are $\{-1,-2,-3\}$-hypomorphic, then they are
$(\leq6)$-hypomorphic and hence $T^{\prime}[X]\sim T[X]$.
\item If $T[X]$ is indecomposable. As $T'[X]$ and $T[X]$ are
$(\leq3)$-hypomorphic, then from Theorem \ref{theo inversion},
$T^{\prime}[X]=T[X]$ or $T^{\prime }[X]=T^{\ast}[X]$. Consider then
the case where $T^{\prime}[X]=T^{\ast}[X]$. If $T[X]$ is self dual,
then $T^{\prime}[X]\sim T[X]$. Suppose now that $T[X]$ is not self
dual. As in addition $T[X]$ is $\{-1,-2,-3\}$-self dual (because
$T^{\prime}[X]$ and $T[X]$ are  \\$\{-1,-2,-3\}$%
-hypomorphic), hence the tournament $T$ belongs to
$O_{3}(I_{n-2,\{-1,-2,-3\}})$ and then to $ \Omega_{n}$; which is
absurd.
\end{itemize}
\end{itemize}
\item \textbf{Fact 5}. If $T \in \Omega_{n}$, then $T' \not\sim T$.

This fact is an immediate consequence of facts 2 and 3.
\end{description}
\hskip 12cm $\square$


\bigskip

We shall now prove Corollary \ref{theo moin3 recon}. Before that let
us observe that Corollary \ref{ajoutth1} is an immediate consequence
of Theorem \ref{Theo moins 3 reconstruction} because each element
$T$ of $\Omega_{n}$ (where $n \geq 9$) has an interval $X$ such that
$T[X]$ is indecomposable and $\mid X \mid \in \{n-1,n-2\}$.


\vskip0.5cm \textbf{Proof of Corollary \ref{theo moin3 recon}.}
Suppose that there exists an integer $n_{0}\geq7$ such that the
indecomposable tournaments with at least $n_{0}$ vertices are
$\{-3\}$-reconstructible and consider a tournament $T$ with $n\geq
n_{0}+2$ vertices. Then, the classes $I_{n-2,\{-3\}}$ and
$I_{n-1,\{-3\}}$ are empty. So, the classes $I_{n-2,\{-1,-2,-3\}}$
and $I_{n-1,\{-2,-3\}}$ are empty. Thus, $\Omega_{n}$ is empty and
then $T \notin \Omega_{n}$. By Theorem \ref{Theo moins 3
reconstruction}, the tournament $T$ is then
$\{-3\}$-reconstructible.

\hskip 11,5cm $\square$

\vskip0.5cm

\center{\emph{AKNOWLEDGEMENT}}\\
We would like to thank the anonymous referee for his insightful
comments.


\bigskip

\begin{thebibliography}{99}                                                                                               %


\bibitem {Bondy}J. A. Bondy, A graph reconstructor's manual, in Surveys in
Combinatorics, (1991), Keendwell \'{e}d., London. Math. Soc. Lecture
Note Ser., Cambridge Univ. Press,  221 - 252.

\bibitem {Bondy et Hemm}J. A. Bondy and R. L. Hemminger, Graph reconstruction,
a survey, J. Graph theory, 1, (1977), 227 - 268.

\bibitem {Boua boud}M. Bouaziz and Y. Boudabbous, La demi-isomorphie et les
tournois fortement connexes finis, C. R. Acad. Sci. Paris S\'{e}rie I Math.
335 (2002), 411-416.

\bibitem {Bouch}H. Bouchaala, Sur la r\'{e}partition des diamants dans un
tournoi, C. R. Acad. Sci. Paris, Serie. I 338 (2004) 109-112.

\bibitem {Bouch et Boud}H. Bouchaala and Y. Boudabbous, La $\{-k\}$%
-autodualit\'{e} des sommes lexicographiques finies de tournois
suivant un $3$-cycle ou un tournoi critique, Ars Combin. 81 (2006).
33-64.





\bibitem {Boud et Bouss}Y. Boudabbous and A. Boussa\"{\i}ri, Reconstruction des
tournois et dualit\'{e}, C.R. Acad. Sci. Paris, t. 320, S\'{e}rie I,
(1995), 397 - 400.

\bibitem {Boud damm ille}Y. Boudabbous, J. Dammak and P. Ille,
Indecomposability and duality of tournaments, Discrete Mathematics 223 (2000) 55-82.





\bibitem {BILT1}A. Boussa\"{\i}ri, P. Ille, G. Lopez and S. Thomass\'{e}, Hypomorphie et inversion locale entre graphes, C. R. Acad. Sci. Paris, S\'{e}rie I 137 (1993), 125-128.

\bibitem {BILT}A. Boussa\"{\i}ri, P. Ille, G. Lopez and S. Thomass\'{e}, The
$C_{3}$-structure of tournaments, Discrete Math., 277 (2004), 29-43.

\bibitem {eur_rosen}A. Ehrenfeucht and G. Rozenberg, Primitivity is hereditary
for $2$-structures, Theoret. Comput. Sci. 3 (70) (1990) 343-358.







\bibitem {Fraiss}R. Fra\"{\i}ss\'{e}, Abritement entre relations et
sp\'{e}cialemnt entre cha\^{\i}nes, Symposi. Math. Instituto
Nazionale di Alta. Mathematica, 5, (1970) 203-251.

\bibitem {Gallai}T. Gallai, Transitiv orientierbare Graphen, Acta. Math. Acad.
Sci. Hungar., 18, (1967), 25-66.





\bibitem{Har et Palm} F. Harary and E. Palmer, On the problem of
reconstructing a tournament from subtournaments, Monatsch. Math. 71,
(1967), 14-23.


\bibitem {Ille}P. Ille, La reconstruction des relations binaires, C. R. Acad.
Sci. Paris, t. 306. S\'{e}rie I, (1988), 635-638.









\bibitem{Lopez 2res} G. Lopez, Deux r\'esultats concernant la
d\'etermination d'une relation par les types d'isomorphie de ses
restrictions, C. R. Acad. Sci. Paris, t. 274, S\'erie A, (1972),
1525 - 1528.


\bibitem{lopez suite} G. Lopez, Sur la d\'{e}termination d'une relation
par les types d'isomorphie de ses restrictions. C. R. Acad. Sci  Paris
(A) t. 275 (1972) 951-953.


\bibitem {Lopez indéf}G. Lopez, L'ind\'{e}formabilit\'{e} des relations et
multirelations binaires, Z. Math. Logik Grundlag. Math., 24, (1978),
303-317.


\bibitem {Lopez+Rauzy I}G. Lopez and C. Rauzy, Reconstruction of binary relations
from their restrictions of cardinality 2, 3, 4 and (n-1), I, Z.
Math. Logik Grundlag. Math., 38, (1992), 27-37.

\bibitem {Lopez+Rauzy II}G. Lopez and C. Rauzy, Reconstruction of binary
relations from their restrictions of cardinality 2, 3, 4 and (n-1),
II, Z. Math. Logik Grundlag. Math., 38, (1992), 157-168.

\bibitem {Moon}J. W. Moon, Topics on Tournaments, Holt, Rinehart and Winston,
New York (1968).

\bibitem {Pouz prop comb}M. Pouzet, Application d'une propri\'{e}t\'{e}
combinatoire des parties d'un ensemble aux groupes et aux relations,
Math. Z. 150 (1976) 117-134.

\bibitem {Pouz rel non reconst}M. Pouzet, Relations non reconstructibles par
leurs restrictions. Journal of Combin. Theory B. vol. 1.26, (1979),
22-34.

\bibitem {Reid et Thom}K. B. Reid and C. Thomassen, Strongly Self-Complementary
and Hereditarily Isomorphic Tournaments, Monatshefte Math. 81,
(1976), 291-304.





\bibitem {Sghiar}M. Sghiar, Private communication.

\bibitem {Stock}P. K. Stockmeyer, The falsity of the reconstruction conjecture
for tournaments, J. Graph Theory 1 (1977), 19-25.

\bibitem {Ulam}S. M. Ulam, A collection of Mathematical Problems,
Intersciences Publishers, New York (1960).
\end{thebibliography}
\end{document}